\documentclass[12pt]{article}

\usepackage[utf8]{inputenc}       
\usepackage[T1]{fontenc}          
\usepackage{lmodern}              
\usepackage{geometry}             
\usepackage{graphicx}             
\usepackage{amsmath, amssymb}     
\usepackage[numbers]{natbib}      
\usepackage{hyperref}             
\usepackage{caption}              
\usepackage{pgfplots}
\usepackage{tikz}
\usepackage{subcaption}

\newcommand{\epsrel}{\varepsilon_{\text{rel}}}
\newcommand{\epsabs}{\varepsilon_{\text{abs}}}

\title{Using a MIP Solver as a PDHG-Based MIP Heuristic}
\author{Edward Rothberg \\
Gurobi Optimization, LLC \\
\texttt{rothberg@gurobi.com}}
\date{September 2, 2026}

\begin{document}

\maketitle
\begin{abstract}
  The PDHG algorithm provides a new capability for solving difficult
  linear programming (LP) problems: the ability to find low-accuracy
  solutions quickly.  While such solutions may not be appropriate for
  all applications of LP, one possible use is to accelerate heuristics
  for finding feasible solutions to Mixed-Integer Programming (MIP)
  problems, where these approximate LP solutions can hopefully
  guide a heuristic towards accurate and high-quality MIP solutions.  While
  a monolithic heuristic that exploits PDHG solutions would be useful,
  we pose a broader question here: could we replace the
  default LP solver in a modern MIP solver with PDHG to accelerate most
  (or all) of its existing heuristics?  We find that doing so does
  have some obvious drawbacks, preventing us from using several
  powerful MIP techniques, but it leaves most others unaffected,
  ultimately resulting in a heuristic that often finds high-quality,
  high-accuracy solutions faster than current state-of-the-art
  strategies.
\end{abstract}

\section{Introduction}

The Primal-Dual Hybrid Gradient (PDHG) algorithm and its recent
enhancements~\cite{pdlp22,pdhg11,hprlp24,cupdlp23,cupdlpx25} have
created intriguing new options for solving linear programming (LP)
problems.  When run on modern, highly parallel computing platforms,
these methods often find low-accuracy solutions to difficult LP
problems much faster than the traditional alternatives: the simplex
method~\cite{simplex51} and interior-point methods~\cite{wright97}.
While solutions with significant constraint violations may not be
appropriate in many contexts, they still contain valuable, high-level
information that could be used to guide other techniques.  In
particular, in a Mixed-Integer Programming (MIP) context, where
heuristics play a large role, they could be used to guide a heuristic
to promising portions of the solution space.  A few researchers have
explored this possibility~\cite{pdhgheuristic25}, with some promising
results.

One lesson from the evolution of MIP solvers over the past few decades
is that finding high-quality solutions requires more than just one
good heuristic.  A modern MIP solver contains dozens of heuristics, as
well as other powerful tools like cutting planes and bound
strengthening that relentlessly work to improve the current
approximation to the global optimum.  These techniques work
in tandem to find good solutions (and good lower bounds), creating a
robust whole from a set of much less robust
parts~\cite{tippingpoint2007}.

This paper therefore takes a different approach to using PDHG to find
high-quality MIP solutions.  Put simply, it poses the question, ``What
would happen if you replaced the LP solver in a state-of-the-art MIP
solver with PDHG?''  At a high level, the answer is that the MIP solver
will be confronted with much lower-accuracy solutions than it was built
to expect.  Some techniques will simply not function after the change;
others may not be as effective.  However, by staying within the framework
of a powerful MIP solver, the parts that remain can benefit
from the interplay between techniques that has proven to be so
effective for MIP solvers in general.

This paper takes a highly quantitative look at this question,
measuring the impact of the swap in relaxation solvers on the various
parts of the MIP solver for a large set of MIP models (the MIPLIB~2017
set~\cite{gleixner2021miplib}).  We find that indeed several powerful
MIP techniques are degraded or cannot be applied when using PDHG
solutions, and this does have a performance impact.  Once we accept
this degradation, though, we find that using low-accuracy solutions
has only a modest impact on the ability of the MIP solver to find
high-quality solutions.  We are left with a powerful tool that can
exploit situations where PDHG is much faster than the alternatives at
solving the LP relaxations.

When we compare the results from our approach to default Gurobi over a
broad set of test models, we find that the former is usually less
effective.  There are several reasons for this, the main one being
that PDHG needs to be substantially faster than the alternatives to
overcome the degradations noted above.  That is unusual on this
test set, probably because most of the test models are not that large.
However, we believe that average performance across a test set is not
the best way to evaluate the effectiveness of a new heuristic; a
heuristic should instead be evaluated based on its ability to expand
the efficient frontier when considering solution quality versus
runtime.  The new heuristic succeeds on that front, finding
high-quality solutions much faster than the alternatives for many
models in the MIPLIB set and also for several customer models.

\section{Background}
\label{sec:Background}

A Mixed-Integer Programming (MIP) model is a problem of the form:
\begin{align*}
\min \quad & c^\top x \\
\text{s.t.} \quad & A x = b \\
& \text{lb} \leq x \leq \text{ub} \\
& \text{some}\ x_j \in \mathbb{Z}
\end{align*}

The most powerful tool for solving such problems is the {\em LP
  relaxation\/}, which is a simplification of the original MIP model
where the integrality restrictions on variables are dropped to create
an LP model.  The solution to that LP model provides a wealth of
information about the MIP model.  One important example is the optimal
relaxation objective value, which provides a bound on the best
possible objective for the associated MIP model.  Another is a set of
{\em reduced costs\/}, which provide lower bounds for each variable
on the degradation of this optimal objective value per unit change in
the value of the variable (from its value in the optimal LP solution).

\subsection{Relaxation Solution Quality}
\label{sec:Quality}

The quality of the information provided by the LP relaxation is only
as good as the quality of the relaxation solution.  Simplex and
interior-point methods provide highly accurate solutions (although
interior-point methods are typically followed by a
crossover~\cite{crossover91} step to achieve such accuracy).  Given a
primal residual vector $r_P = A x - b$ that captures constraint
violations, the default termination tolerance for simplex and
crossover is typically $\|r_P\|_{\infty} \leq \epsabs$ (i.e., no
constraint is violated by more than $\epsabs$), with
$\epsabs = 10^{-6}$ being traditional.

Termination tolerances are typically much looser for PDHG, with
$\|r_P\|_2 \leq \epsrel (1 + \|b\|_2)$ being typical.  A default value
of $\epsrel = 10^{-4}$ is common in academic codes, although
commercial implementations appear to have settled on
$\epsrel = 10^{-6}$.  While $10^{-6}$ values for $\epsabs$ and
$\epsrel$ may appear similar, the PDHG tolerance is dramatically
looser; we have found that PDHG produces maximum absolute constraint
violations that are nearly six orders of magnitude larger on
average~\cite{pdhghybrid26}.  To make this more concrete in our
context, a constraint $x_1 + x_2 \leq 1$ that enforces a mutual
exclusion between a pair of binary variables will always be satisfied
to within $10^{-6}$ in a simplex or interior-point solution of the
relaxation, while a PDHG solution can completely ignore a substantial
number of such constraints while still satisfying its termination
criteria.  It will come as no surprise that this difference in
accuracy will be important when substituting PDHG for the default LP
solver in a MIP solver.

One other difference between PDHG solutions and the solutions the MIP
solver is normally accustomed to is that PDHG solutions are not {\em
  basic\/}.  Given an $m \times n$ constraint matrix $A$, a basic
solution will have at least $n-m$ variables at a bound.  Basic
solutions are a large part of the reason why simplex solutions are so
much more accurate, and they provide other benefits that will be
discussed later.

\section{MIP Solver Structure}
\label{sec:Structure}

Before tackling the question of how substituting PDHG for the default
LP solvers will affect the overall MIP solution process, we need to
provide a lot more information about the various pieces of a MIP
solver.  We give a brief overview of the most important pieces for our
purposes here, but we refer the reader to a more comprehensive
discussion (e.g.,~\cite{wolsey2020}) for additional information.

Solving a MIP is a process of continuous improvement.  At any point in
the process, the solver maintains a lower bound on the optimal
objective value (assuming minimization), which typically comes from
solving relaxations, and an upper bound, which comes from the best
feasible solution that has been found so far (referred to as
the {\em incumbent}).

The two main tools available for improving the lower bound are cutting
planes~\cite{wolsey2020} and branching~\cite{land1960}.  Cutting
planes are linear constraints that cut off a portion of the feasible
region for the LP relaxation without affecting the set of feasible
solutions to the MIP model.  Modern MIP solvers include many varieties
of cutting planes, and they can be quite powerful at improving the
bound.

Branching splits the feasible region into two or more pieces by adding
constraints.  The typical approach for MIP is to identify
an integer variable $x$ that takes a fractional value $x^*$ in the
current LP relaxation and split that MIP problem into two by
adding $x \leq \lfloor x^* \rfloor$ in one child and
$x \geq \lceil x^* \rceil$ in the other.  By excluding the relaxation
solution, the hope is that the relaxation objective values for the two
children will both be larger than that of the parent, thus increasing
the overall lower bound.  This process continues recursively, dividing
the original solution space into a {\em branch-and-bound tree\/} of
disjoint subspaces.

While branching will eventually produce feasible solutions for the MIP
model (assuming it is feasible), one lesson learned from many decades
of work on MIP solvers is that it is best to find good feasible
solutions as early in the solution process as possible.  A modern MIP
solver employs a variety of heuristics to find such solutions.

Another lesson learned is that it is best to make as much progress
as possible toward closing the gap between the lower and upper bounds
before branching starts.  Thus, a modern MIP solver will expend a
substantial amount of work at the {\em root node\/} of the
branch-and-bound tree.

\begin{figure}[htpb]
\centering
\begin{tikzpicture}[
    node distance=1.6cm,
    block/.style={
        rectangle,
        draw,
        rounded corners,
        minimum width=3.0cm,
        minimum height=1.0cm,
        align=center
    },
    arrow/.style={->, thick}
]

\path[use as bounding box] (-2.2,-5.8) rectangle (2.2,0.8);

\node[block, fill=blue!12, draw=blue!60!black] (pre) {Preprocessing};
\node[block, fill=blue!12, draw=blue!50!black, below of=pre] (relax) {Relaxation solve};
\node[block, fill=blue!15, draw=blue!70!black, below of=relax] (cuts) {Root cut loop};
\node[block, fill=blue!12, draw=blue!60!black, below of=cuts] (tree) {Tree exploration};

\draw[arrow] (pre) -- (relax);
\draw[arrow] (relax) -- (cuts);
\draw[arrow] (cuts) -- (tree);

\draw[arrow] (cuts) edge[loop right] (cuts);
\draw[arrow] (tree) edge[loop right] (tree);

\end{tikzpicture}
\caption{MIP solver phases}
\label{fig:mip_phases}
\end{figure}
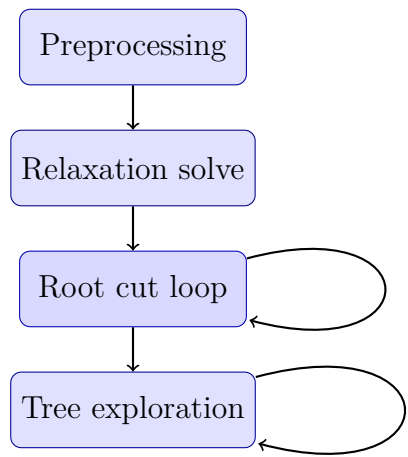

Let us now take a closer look at the major phases of computation in a
MIP solver, illustrated in Figure~\ref{fig:mip_phases}.

\subsection*{Preprocessing}

The first phase in a MIP solve is {\em preprocessing\/}, where the
original MIP model is transformed into a new model using a variety of
reductions that create a smaller and {\em tighter\/} model that has
the same optimal objective value as the original model
(see~\cite{presolve2020} for examples of some of the more involved
presolve reductions).  MIP heuristics are often run during this phase,
on both the original and presolved models, with the goal of
obtaining a feasible solution very early in the solution process.
These heuristics cannot take advantage of relaxation solution
information, so their quality is typically (but not always) quite
poor.

\subsection*{Solving the relaxation}

The second phase is the {\em relaxation solve\/}, where a solution is
computed for the LP relaxation of the presolved model.  As noted
earlier, the solution to this model typically provides essential
guidance for the steps that follow, including for the heuristics that
aim to find high-quality solutions to the MIP model.

The root relaxation is typically solved with either the simplex method
or an interior-point method (or both running concurrently).  This
choice is automated in Gurobi by default.  For a variety of reasons,
including accuracy, reoptimization speed, and conciseness of solution
representation, Gurobi always computes basic optimal solutions to LP
relaxations.  An interior-point solve is therefore always followed
by a crossover step.

\subsection*{Root cut loop}

Once the relaxation has been solved, the next phase performs the {\em
  root cut loop\/}, generating cutting planes to tighten the
relaxation.  Gurobi 13.0 contains 26 different cutting plane types,
and after several root cut loop iterations these typically
substantially reduce the gap between the root relaxation objective
value and the optimal solution objective value.

The root cut loop is made up of a number of steps that will be
discussed in more detail shortly.  At this point, we just want to
reiterate that modern MIP solvers put a lot of effort into trying to
close as much of the optimality gap as possible at the root node.

\subsection*{Tree exploration}

The final phase in the solution of a MIP model is exploring the
branch-and-bound tree.  Most of the techniques applied at the root
node (including heuristics, cutting planes, etc.) are also applied at
the nodes of the branch-and-bound tree, but for various reasons
they typically do not have the same impact when applied at the nodes.

\section{Experimental Methodology}

This paper reports on extensive computational experiments.  We measure
the effectiveness of the various phases of the MIP computation at
finding good feasible solutions, the impact of disabling various MIP
features, and the ability of our PDHG-based approach running on a GPU
to find high-quality MIP solutions faster than the alternatives.  This
section briefly discusses our computational testing methodology.

\subsection{Software and hardware}

All of our computational tests use a modified version
of the Gurobi Optimizer (version 13).  One set of modifications adds
instrumentation that enables us to gather and report relevant
statistics.  Another allows the code to use PDHG (without crossover,
using $\epsrel = 10^{-6}$) to solve LP relaxations.  Finally, we add
parameters to allow certain features to be shut off, to support
ablation testing.

When comparing our PDHG-based approach against the
alternative, we always use the concurrent LP solver to solve the root
relaxation (by setting the Gurobi {\em Method\/} parameter to value
4).  This choice ensures that the automatic scheme does not
miss a simple opportunity to improve performance.
We refer to this configuration as {\em baseline Gurobi\/}.

We use two machines for our tests.  In almost all cases, CPU testing
is performed on an AMD EPYC 4364P system with 8-core processors and
128~GB of DDR5 memory.  We use a 10,000-second time limit (wall-clock time)
for our baseline runs, and a 50,000-second time limit for some PDHG tests
(since PDHG runs particularly slowly on this system).

We perform our GPU tests and some of our CPU tests on a system with an
AMD EPYC 9575F CPU and an NVIDIA RTX PRO 6000 GPU.  The CPU has 64
cores, the machine has 768~GB of system memory, and the GPU has 96~GB
of fast GDDR7 memory.  For reference, the GPU is roughly 40 times faster
than the EPYC 4364P CPU when performing PDHG iterations for large
LP models.

\subsection{Test sets}

As noted earlier, we use the $1,065$ models of the MIPLIB 2017 test
set for our computational testing.  This collection has significant
limitations for our purposes, so many of our tests only consider the
$87$ models in this set where the MIP root node finds a feasible
solution and requires between $100$ and $10,000$ seconds to complete
on our EPYC 4364P system using baseline Gurobi.  This set is further
refined in one test to a set of $68$ models where the root node
requires no more than $50,000$ seconds for our PDHG-based heuristic
approach, and all methods considered find at least one feasible
solution.  We also report detailed results for a set of 12 models
where our approach performs particularly well, as well as results for
a few anonymized customer models.

\section{Why Focus on the Root Node?}

The focus of this paper is heuristics, so let us now sharpen our focus
by taking a quantitative look at which phases of the MIP solution
process are most effective at finding high-quality feasible
solutions.  Specifically, we now measure the quality of the best
feasible solution available at the end of each MIP phase for the
MIPLIB 2017 test set, running baseline Gurobi on our EPYC 4364P systems
with a time limit of $10,000$ seconds.  We start by reporting the
number of models that have feasible solutions at the end of each
phase.  For those models, we then capture the best objective found at
the end of that phase ($z_{\text{phase}}$), and the best objective
found in 10,000 seconds ($z_{\text{best}}$).  We compute the ratio
$(z_{\text{phase}} - z_{\text{best}}) / \max(1.0, |z_{\text{best}}|)$
for each model, and then take the geometric mean of this value over
all models where feasible solutions were found in that phase.  Gaps
smaller than $10^{-4}$ are rounded up to $10^{-4}$.  Models are
grouped into those where at least $10$ tree nodes are explored before
the MIP solve terminates, those where at least $1,000$ nodes are
explored, and those where at least $100,000$ nodes are explored.

Table~\ref{tab:heur_improvement} shows the results of this test.
\begin{table}[htbp]
\centering
\begin{tabular}{c|c|c|c|c|c|c|c}
   & \multicolumn{7}{|c}{Solution status after each phase} \\
  Minimum nodes & \multicolumn{2}{c|}{Preprocessing} & \multicolumn{2}{c|}{One cut pass} & \multicolumn{2}{c|}{Root} & MIP search \\
  explored & Sols & Gap & Sols & Gap & Sols & Gap & Sols \\ \hline
  10 & 431 & 262\% & 556 & 20\% & 618 & 1.8\% & 689 \\
  1,000 & 379 & 254\% & 486 & 23\% & 540 & 2.4\% & 602 \\
  100,000 & 169 & 247\% & 209 & 25\% & 231 & 5.2\% & 252\\
\end{tabular}
\caption{Solution statistics (number of models where solutions are found and geometric mean gap of those solutions) after various MIP phases.}
\label{tab:heur_improvement}
\end{table}
The table includes results for the phases defined so far plus two more
sets of results.  The first additional results are for ``One cut
pass'', which captures results after one pass of the root cut loop,
using all of the heuristics that would be applied before the first
reoptimization.  The second set
of additional results are for ``MIP search'', which shows the number
of models where solutions are found before termination (where termination
is due to the 10,000-second time limit for 313 models and due to
successful termination of the MIP search for the rest).  Note that the
gap at this point will always be 0, since this is the source of our
best known incumbent value.

As the results show, the root node heuristics do the bulk of the work
of finding high-quality solutions.  For models that explore at least
10 MIP nodes, the mean gap for the solution available at the end of
the root is less than $2\%$.  Even for the models that explore at
least $100,000$ MIP nodes, the gap is just over $5\%$.  Gaps after
a single root cut pass are roughly 10 times worse, and
gaps after preprocessing are another roughly 10 times worse.

Because of the results of this test, we set a goal of using PDHG
solutions to perform the MIP solution process through the end of the
root cut loop.

\section{Root Cut Loop}
\label{sec:RootCutLoop}

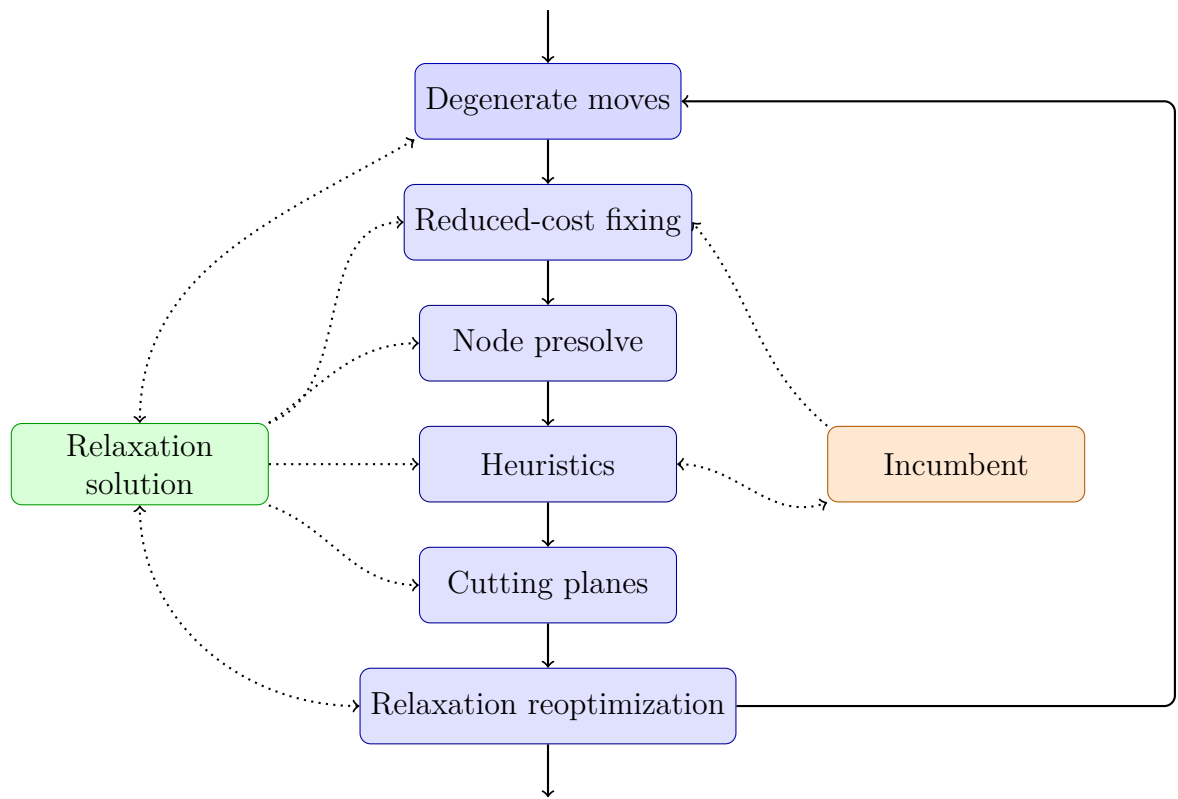
\begin{figure}[htpb]
\centering
\begin{tikzpicture}[
node distance=1.6cm,
block/.style={
rectangle,
draw,
rounded corners,
minimum width=3.4cm,
minimum height=1.0cm,
align=center
},
arrow/.style={->, thick},
bidir/.style={<->, thick},
sidearrow/.style={->, thick, dotted},
sidebidir/.style={<->, thick, dotted}
]

\path[use as bounding box] (-7.9,-9.5) rectangle (7.9,1.4);

\node[block, fill=blue!15, draw=blue!70!black] (deg) {Degenerate moves};
\node[block, fill=blue!12, draw=blue!60!black, below of=deg] (rcf) {Reduced-cost fixing};
\node[block, fill=blue!12, draw=blue!50!black, below of=rcf] (presolve) {Node presolve};
\node[block, fill=blue!12, draw=blue!50!black, below of=presolve] (heur) {Heuristics};
\node[block, fill=blue!12, draw=blue!60!black, below of=heur] (cuts) {Cutting planes};
\node[block, fill=blue!12, draw=blue!60!black, below of=cuts] (resolve) {Relaxation reoptimization};

\node[block, fill=green!15, draw=green!60!black] (relaxsol) at (-5.4,-4.8) {Relaxation\\solution};
\node[block, fill=orange!18, draw=orange!70!black] (incumbent) at (5.4,-4.8) {Incumbent};

\draw[arrow] (deg.north) ++(0,0.7) -- (deg.north);

\draw[arrow] (deg) -- (rcf);
\draw[arrow] (rcf) -- (presolve);
\draw[arrow] (presolve) -- (heur);
\draw[arrow] (heur) -- (cuts);
\draw[arrow] (cuts) -- (resolve);

\draw[arrow] (resolve.south) -- ++(0,-0.7);

\draw[arrow, rounded corners] (resolve.east) -- ++(5.8,0) |- (deg.east);

\draw[sidebidir] (resolve.west) to[out=180,in=-90] (relaxsol.south);
\draw[sidebidir] (deg.south west) to[out=210,in=90] (relaxsol.north);
\draw[sidearrow] (relaxsol.north east) to[out=20,in=180] (rcf.west);
\draw[sidearrow] (relaxsol.north east) to[out=35,in=180] (presolve.west);
\draw[sidearrow] (relaxsol.east) to[out=0,in=180] (heur.west);
\draw[sidearrow] (relaxsol.south east) to[out=-20,in=180] (cuts.west);

\draw[sidearrow] (incumbent.north west) to[out=140,in=-50] (rcf.east);
\draw[sidebidir] (heur.east) to[out=0,in=200] (incumbent.south west);

\end{tikzpicture}
\caption{MIP root node steps.}
\label{fig:mip_root}
\end{figure}

To understand the impact of using PDHG to solve relaxations, we now
examine the question of how these solutions will affect the
various steps of the root cut loop.  We now briefly describe what each
of these steps does and how it uses the relaxation.  Recall that while
adding cutting planes to the relaxation is the main
motivation for the root cut loop, this phase includes a
number of other steps, illustrated in Figure~\ref{fig:mip_root}.
The current relaxation and incumbent solutions are represented as
boxes on the left and right of the figure.  The various steps in the
root cut loop are shown in the middle, with directed edges indicating
whether they produce or consume relaxation solutions or incumbents.

\subsection*{Degenerate moves}

An LP problem typically has multiple optimal solutions, often with
significantly different characteristics (the most interesting in our
context being the number of integer variables from the MIP that take
fractional values in the relaxation solution).  The {\em degenerate
  moves\/} step explores the set of optimal solutions, looking for
relaxation solutions with fewer integer-infeasible variables.  Such
solutions often make it much easier for heuristics to find feasible
MIP solutions.

\subsection*{Reduced-cost fixing}

Recall that {\em reduced costs\/} provide lower bounds on the rate of
change in the objective value relative to a change in the value of a
variable.  When used in conjunction with the current incumbent
objective, reduced costs can be used to establish a value for each
variable at which the relaxation solution objective would be no better
than that of the current incumbent, and thus not of interest.
This can lead to a new bound on the value of that variable.

\subsection*{Node presolve}

Given new information derived during the root cut loop (e.g., new
bounds on variables or new cutting planes), it is often possible to
derive additional variable bounds through simple propagation of linear
constraints.  To give a simple example, given a constraint $x \leq 10 b$
(with $x$ non-negative and $b$ binary), if it can be shown (e.g.,
through reduced-cost fixing) that $b=0$, it follows that $x = 0$.

\subsection*{Heuristics}

Gurobi invokes a number of primal feasibility heuristics in each root
cut pass.  Gurobi 13.0 implements 52 different heuristics; the subset
that are invoked in each cut pass will vary, depending on
characteristics of the problem and information learned from previous
attempts to apply these heuristics.  Heuristics fall into two general
categories: construction heuristics, which attempt to construct a new
incumbent, and improvement heuristics, which attempt to improve the
current incumbent.  In all cases, heuristics rely heavily on the
current relaxation solution to guide the search for a new incumbent.

Most of the 52 Gurobi heuristics are {\em sub-MIP\/} heuristics,
meaning that they fix variables in the original MIP model and solve a
smaller MIP model recursively.

While each of these heuristics plays a role, the most effective by far
is RINS~\cite{rins2005}.  RINS is an improvement heuristic that uses
the current incumbent solution and the current relaxation to identify
variables that are good candidates for fixing, leaving a sub-MIP that
typically captures a productive neighborhood of the current incumbent.
Note that RINS never works alone.  As an improvement heuristic, it
requires a feasible solution to start from, and it often improves on
feasible solutions found by other heuristics.

\subsection*{Relaxation reoptimization}

A single pass through the root cut loop can change the MIP
problem dramatically.  Cutting planes can add constraints,
reduced-cost fixing and node presolve can add new bounds, etc.  At the
end of each root cut pass, we need to reoptimize the LP relaxation to
incorporate this new information and give the next pass a new
relaxation solution to work from.  This reoptimization is typically
performed with the dual simplex method, since the changes made to
the model preserve dual feasibility, and in particular preserve dual
feasibility of the current basis.  This is advantageous for dual
simplex.

\subsection*{Exploiting variability with multiple concurrent root solves}

One technique used in the root cut loop that is not reflected in
Figure~\ref{fig:mip_root} is the use of multiple threads to exploit
performance variability.  Several steps in the root cut loop
(heuristics in particular, and also cutting planes to some extent) can
produce significantly different results after small perturbations of
their inputs.  It is therefore beneficial to start multiple threads,
each working on a perturbed version of the root relaxation in
parallel.  These threads share any new feasible solutions and cutting
planes that they find along the way.

\subsection*{Interplay between algorithms}

It should be apparent from the discussion so far that these methods
never work alone.  The output of one feeds into the input of another,
sometimes within one root cut loop iteration and other times across
iterations.  For example, a new incumbent can lead to new bounds from
reduced-cost fixing and node presolve, which can change the relaxation
solution.  The new relaxation can send construction heuristics down a
new path, and a new incumbent can similarly impact the construction
heuristics.  This cooperation is at the heart of what makes modern MIP
solvers so effective.

\subsection*{Runtime costs}

Let us add a brief note about the relative costs of the various
components described above.  We have found through extensive profiling
that reoptimizing the relaxation is almost always the most expensive
step.  Part of this is by design, since the runtimes of the other
steps (particularly heuristics and cutting planes) can be dialed up or
down by simply stopping them before they exhaust all possibilities.
Extensive testing has produced little evidence that letting these run
for longer would lead to significant improvements.

The fact that solving LPs is the dominant cost is good news for our
proposed PDHG-based heuristic; if PDHG is able to solve LP relaxations
much faster, that creates significant scope for speeding up the overall
root cut loop.

\section{What Do We Lose with PDHG?}

Recall that our goal is to exploit the ability of PDHG to find
low-accuracy solutions to LP relaxations quickly,
thereby accelerating the root cut
loop and finding high-quality MIP solutions more quickly.  We now consider the
question of what is lost by using these low-accuracy, non-basic
relaxation solutions.  We consider this question in three parts: (i)
which steps in the root cut loop become less effective, (ii) which
potentially get slower, and (iii) which are simply no longer possible.
We try to quantify the impact of each.

One point we should make before proceeding is that while the lack of
accuracy in PDHG solutions would appear to make them unsuitable for
computing lower bound information (relaxation objective bounds,
reduced-cost based bounds, etc.), such solutions can in fact be used
for this purpose in many cases.  Specifically, valid bounds can always
be obtained when all variables in the model have finite lower and
upper bounds (since this effectively makes any dual solution feasible).
This condition is met for many important classes of
models (e.g., pure binary models).  Note that the resulting bound
information may not be as strong as it would be for a true optimal
solution.

\subsection{Which Steps Become Less Effective?}

Let us first consider which of the steps performed as part of the root
cut loop could become less effective when using a PDHG solution.  The
two that probably seem the most worrisome are the two most important
steps: cutting plane separation and heuristics.  Most of the
heuristics and all of the cutting plane separators exploit the current
relaxation solution, and a less accurate solution could clearly have
an impact.  At the same time, these steps are not necessarily that
sensitive to the accuracy of the solution.  Heuristics in particular
are guided by the relaxation solution, but the guidance is often quite
coarse.  We will come back to this point, but for now let us note that
this is an area of concern.

\subsection{Which Steps Get Slower?}

The main step in the root cut loop that has the potential to get
slower from switching to PDHG is the LP reoptimization, which is
ordinarily performed using dual simplex.  Simplex is renowned for its
ability to reoptimize an LP quickly after small changes.  PDHG can
also be warm-started from a previous solution, but the cost is less
well studied.  Since LP reoptimization is the dominant cost of the
root cut loop, slower reoptimization could have a significant impact
on runtime.

We have tried to quantify this by looking at the geometric mean of
the ratio of the time required for the reoptimization at the end of
each root cut pass relative to the time required to solve the original LP
relaxation from scratch (using $10^{-2}$ when the computed ratio is
smaller).
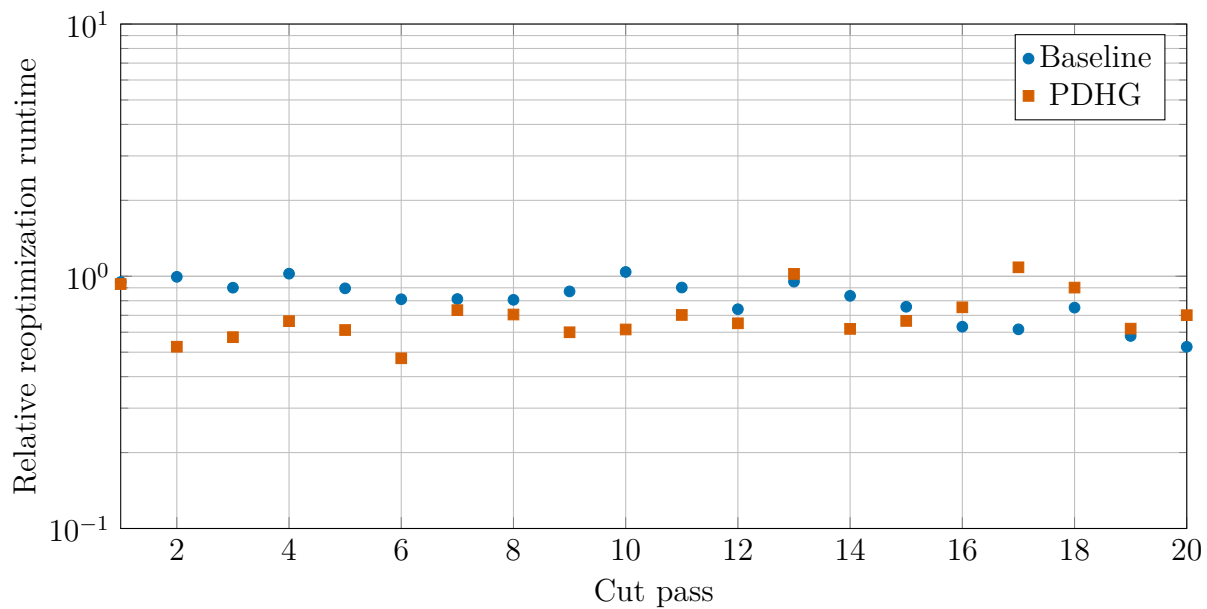
\begin{figure}[htbp]
\centering
\begin{tikzpicture}
\begin{semilogyaxis}[
    width=0.95\linewidth,
    height=0.5\linewidth,
    grid=both,
    legend style={at={(0.98,0.98)},anchor=north east},
    mark size=1pt,
    xmin=1,
    xmax=20,
    ymin=1e-1,
    ymax=1e+1,
    xlabel={Cut pass},
    ylabel={Relative reoptimization runtime},
    legend image post style={mark size=2pt},
    ]

\definecolor{OIblue}{RGB}{0,114,178}
\definecolor{OIorange}{RGB}{230,159,0}
\definecolor{OIteal}{RGB}{0,158,115}
\definecolor{OIvermillion}{RGB}{213,94,0}

\addplot[
    only marks,
    mark=*,
    mark size=2pt,
    color=OIblue,
] table {
x y
1 0.948878
2 0.994923
3 0.900907
4 1.02406
5 0.89573
6 0.810659
7 0.812301
8 0.806529
9 0.870983
10 1.0399
11 0.90236
12 0.740272
13 0.95218
14 0.836459
15 0.757002
16 0.631115
17 0.616168
18 0.75124
19 0.579881
20 0.52518
};
\addlegendentry{Baseline}

\addplot[
    only marks,
    mark=square*,
    mark size=2pt,
    color=OIvermillion,
] table {
x y
1 0.931552
2 0.525595
3 0.5736
4 0.663965
5 0.612259
6 0.473011
7 0.733897
8 0.705352
9 0.599791
10 0.615297
11 0.702519
12 0.651024
13 1.02187
14 0.618623
15 0.665382
16 0.753326
17 1.0852
18 0.901764
19 0.619976
20 0.701155
};
\addlegendentry{PDHG}

\end{semilogyaxis}
\end{tikzpicture}
\caption{Relaxation reoptimization time relative to root relaxation solve time.}
\label{fig:root_solve_time}
\end{figure} We used our 87-model MIPLIB subset for this
test.  Results for this subset are shown in
Figure~\ref{fig:root_solve_time}.

Recall that we have forced Gurobi to use the concurrent LP solver to
solve the original root relaxation in our tests.  In contrast, Gurobi
will nearly always use dual simplex to reoptimize an LP after it has
been modified.  Thus, the {\em Baseline\/} ratios shown in the figure
compare dual simplex in the numerator against the fastest of primal
simplex, dual simplex, and parallel interior-point in the denominator.

The {\em Baseline} line in the figure shows that reoptimization times
are quite consistent as root cut iterations proceed under baseline
settings, and are slightly less expensive on average than the solve on
the original relaxation.  There are two dynamics at work here.  First,
as iterations proceed, more cutting planes are added, making the
problem larger and thus more expensive to solve.  On the other hand,
the number of cutting planes added in each cut round generally
decreases, leading to smaller changes in the model and thus easier
reoptimization tasks.  These appear to balance each other out, leading
to consistent overall costs.

The {\em PDHG} line in the figure shows the cost of reoptimizing with
PDHG relative to the cost of solving the initial relaxation with PDHG.
The two lines in the figure suggest that the relative reoptimization
costs for PDHG and dual simplex are quite comparable.  In other words,
PDHG appears to also be quite effective at reoptimizing LP
relaxations.

\subsection{Which Steps Are No Longer Possible?}

Now let us consider the question of which root cut loop steps are no
longer possible if we use PDHG to solve relaxations.  The problematic
properties of PDHG are lack of accuracy, lack of a basis, and the
fact that PDHG is a heavy consumer of computing resources (CPU cores
or GPUs).  We discuss the affected steps first and then try to
quantify the effect of not being able to use them.

If we consider the steps in Figure~\ref{fig:mip_root}, we find that
the degenerate move step is the first to be affected.  This step is
usually fast, but its speed depends on the ability to formulate a
smaller LP whose solutions are all optimal for the original LP
relaxation.  This smaller LP can typically be solved quickly.  In the
context of PDHG, previous work~\cite{pdhgcorner25} has found that,
without an unambiguously optimal solution to the LP relaxation to
start from, formulating a smaller problem that maintains optimality
becomes difficult, which makes it expensive to explore the space of
alternate optimal solutions.  We therefore view this step as being impractical
when PDHG is the LP solver.

Most of the cutting plane routines generate violated cuts from an
arbitrary solution vector.  That vector does not need to be optimal
for the LP relaxation (or even feasible); the resulting cutting planes
will still be valid.  There are a few cutting planes, particularly
Gomory cuts, whose separation procedures require a basis.  Without the
information contained in a basis, they simply cannot function and have
to be disabled.

The last root cut step that is not possible with PDHG is the use of
multiple threads to exploit performance variability.  One
often-overlooked advantage of using simplex as the LP solver in MIP is
the fact that it only uses one core, which opens up opportunities to
exploit parallelism.  By contrast, PDHG requires substantial
parallelism to get good performance, using multiple cores on a CPU or
an entire GPU.  If we used PDHG to solve relaxations, launching
several such solves in parallel would exhaust available resources
quickly.  We therefore only launch a single root cut loop.

It may appear that reduced-cost fixing is not possible (or advisable)
with the low-accuracy solutions produced by PDHG.  While basing
variable fixing decisions on noisy reduced-cost information can
certainly lead to invalid fixings, our intent is to
use the root cut loop to generate heuristic solutions, and invalid
fixings will not lead to invalid solutions.  Also, as mentioned
earlier, valid bound information can be obtained from low-accuracy
solutions for many common classes of models.

\subsection{Other Considerations}

We next discuss a few other relevant issues that arise when using
PDHG as the engine for solving relaxations within a MIP solver.

We mentioned that most of our heuristics are sub-MIP heuristics,
meaning that they (repeatedly) solve smaller MIP models.  A natural
question is how to solve the relaxations that arise in these smaller
MIPs.  In this paper, we revert to using the default LP relaxation
algorithm.  These sub-MIP models are constructed to be much smaller
than the original MIP model, so they are much more likely to be
efficiently handled by the default approach.  This strategy could be
problematic for truly enormous MIP models; selectively using
PDHG for some sub-MIP models is a possible remedy.
We have not explored that question here.

Another important issue for PDHG is model scaling.  Adding cutting
planes to the model adds constraints to the relaxation, which changes
the conditioning of the constraint matrix.  Scaling is important for
simplex and absolutely crucial for PDHG.  When using simplex, the
existing matrix scaling is simply extended to accommodate new rows,
which has been found to avoid major conditioning issues.  However, we
have found that this approach is not sufficient for PDHG, and have
chosen to rescale the entire matrix after each cut round.  The cost of
this rescaling is roughly equivalent to the cost of a few dozen PDHG
iterations, which is almost always small in comparison to the number
of iterations required for reoptimization.

Another issue that must be mentioned for PDHG is the termination
criterion.  We use the Gurobi default of $\epsrel = 10^{-6}$, but our
goal in this paper is to consider the impact of low-accuracy
solutions, so we could consider using looser tolerances.  We found
that relaxing the tolerance (e.g., to $\epsrel = 10^{-4}$) quickly led
to a problem.  Cutting planes are new constraints that are violated by
the current relaxation solution.  Using infinity-norms for
termination tolerances, the simplex solver will always perform
iterations to resolve a violation that exceeds the tolerance in any
constraint.  For the 2-norm measure used in PDHG, we found that the
violation introduced by a set of violated cutting planes was often
small enough that the overall violation was still within the desired
tolerances.  This often led the cut loop to be terminated early, which
harmed its ability to continue to find improving solutions.

One potential concern when using low-accuracy LP solutions is that
they will produce low-accuracy MIP solutions.  This was not an issue
here, for two reasons.  The first is that our most effective
heuristics use a sub-MIP approach, as noted earlier.  This means that
the default LP solver is used to compute the ultimate MIP solution.
The second is that Gurobi by default is quite fastidious about
checking the quality of heuristic solutions.  When violations are
found, the solution is simply discarded.  That happens more often when
using low-accuracy LP solutions, which does represent some amount of
wasted work.

\subsection{Comparison}

Now that we have listed the features that we have had to disable to
accommodate PDHG solutions in the MIP solver, let us quantify the
performance impact of removing these now-missing features.  We do that
by first disabling them, one by one, from the baseline Gurobi MIP
solver.  Table~\ref{tab:remove_features1} shows results for this
ablation study, looking at the gap between the solution obtained after
the root node completes for the various degraded versions we consider
and the best solution found in 10,000 seconds.  For this test, we
again use our 87-model subset run on our EPYC 4364P systems.
\begin{table}[htbp]
\centering
\begin{tabular}{c|c|c}
  Additional feature disabled & Models with incumbents & Mean gap to best \\ \hline
  None & 87 & 1.1\% \\
  Multiple root threads & 82 & 3.5\% \\
  Degenerate moves & 83 & 3.6\% \\
  Gomory cuts & 82 & 3.8\% \\
\end{tabular}
\caption{Impact of removing features (over 87 models).}
\label{tab:remove_features1}
\end{table}
The first line in the table shows results for our baseline Gurobi run;
each successive line shows the results obtained by shutting off one
additional feature.  It is apparent that removing these features has
an impact; the number of models where incumbents are found drops from
87 to 82 or 83 (but we should note that a model was only included in
the 87-model test set if baseline Gurobi found an incumbent, so this
decrease may be partially due to variability between runs).  The mean
gap to the best incumbent increases from 1.1\% to 3.8\%.

Now let us add PDHG into the mix.  Table~\ref{tab:remove_features2}
presents similar data, comparing (i) baseline Gurobi, (ii) Gurobi with
the impacted features disabled (labeled {\em Degraded baseline Gurobi\/}),
and (iii) Gurobi using PDHG (with the same impacted features disabled).
\begin{table}[htbp]
\centering
\begin{tabular}{c|c}
  Variant & Mean gap to best \\ \hline
  Baseline Gurobi & 0.8\% \\
  Degraded baseline Gurobi & 3.3\% \\
  PDHG & 4.9\% \\
\end{tabular}
\caption{Impact of removing the three features that had to be disabled for PDHG (over 68 models for which all considered approaches complete the root node within the time limit).}
\label{tab:remove_features2}
\end{table}
For the three-way comparison in Table~\ref{tab:remove_features2},
we further restrict the 87-instance subset to the 68 instances on which
all three variants complete root processing within their respective
time limits and produce at least one incumbent.
We ran this test on our EPYC 4364P systems with the same
10,000-second time limit, except for the PDHG run where we used a
50,000-second time limit (since this machine gives poor performance
for PDHG).

The intent of this test is to be able to focus on the question of how
much impact solution accuracy has on the ability of the various Gurobi
techniques to find high-quality solutions at the root.  The only
substantive difference between the {\em degraded Gurobi\/} and {\em
  PDHG\/} lines in this table is the relaxation solutions that they
use.  Our conclusion from the results in the table is that the
degradation that comes from lower accuracy is quite modest.

\subsection{Summing Up}

To sum up, using PDHG to solve relaxations presents some challenges
but also some significant opportunities.  Once you accept the costs of
shutting off features that are impractical or unavailable when using
PDHG solutions, we find that the quality of the MIP solutions you get
from these solutions is quite similar to what you would get from basic
solutions, despite the dramatic difference in accuracy.  That appears
to provide a promising base on which to build a PDHG-based heuristic.

\section{Computational Results}

We will now stop teasing apart the small differences between baseline
Gurobi and our proposed PDHG-based approach and instead evaluate
results for the proposed approach as a whole.  In other words, we now
look at the question of whether the PDHG-based approach produces
better heuristic solutions in less time.

Aggregating results over the whole MIPLIB set, or even over our
87-model subset, unfortunately does not show an overall advantage for
PDHG.  For PDHG to win, it would have to solve relaxations
significantly faster than Gurobi's baseline algorithm.  Most of the
models in our 87-model subset only spend a few hundred seconds at the
root, and that time is typically split among many cut loop iterations.
The relaxations being solved can hardly be considered difficult LP
models, so they fall outside the sweet spot for PDHG.  Given this, and
given what was said earlier about evaluating heuristics in terms of
their ability to expand the efficient frontier of what is possible, we
present results for several specific instances instead.

\begin{figure}[htbp]
\centering

\begin{subfigure}{0.45\textwidth}
\centering
\begin{tikzpicture}
\begin{axis}[
    grid=both,
    legend style={at={(0.98,0.98)},anchor=north east},
    mark size=1pt,
    xlabel={Runtime (s)},
    ylabel={Incumbent objective},
    ylabel style={yshift=5pt},
    legend image post style={mark size=2pt},
    title={neos-3322547-alsek}
    ]

\definecolor{OIblue}{RGB}{0,114,178}
\definecolor{OIvermillion}{RGB}{213,94,0}

\addplot[
    const plot mark left,
    mark=*,
    mark size=2pt,
    color=OIblue,
] table {
x y
0.07 692
8.9 507
8.94 504
8.97 503
9.04 502
9.6 444
9.63 443
9.7 441
9.91 440
9.99 439
10.99 437
11.65 433
11.89 431
11.95 430
13.63 419
14.57 418
15.68 417
16.61 416
19.02 415
19.38 412
27.85 410
47.61 407
139.28 405
};
\addlegendentry{Baseline}

\addplot[
    const plot mark left,
    mark=square*,
    mark size=2pt,
    color=OIvermillion,
] table {
x y
0.08 692
2.22 460
2.63 413
2.68 403
4.82 402
95.75 401
};
\addlegendentry{PDHG}
\end{axis}
\end{tikzpicture}
\end{subfigure}
\hfill
\begin{subfigure}{0.45\textwidth}
\centering
\begin{tikzpicture}
\begin{axis}[
    grid=both,
    legend style={at={(0.98,0.98)},anchor=north east},
    mark size=1pt,
    xlabel={Runtime (s)},
    ylabel style={yshift=5pt},
    legend image post style={mark size=2pt},
    title={neos-5114902-kasavu},
    ]

\definecolor{OIblue}{RGB}{0,114,178}
\definecolor{OIvermillion}{RGB}{213,94,0}

\addplot[
    const plot mark left,
    mark=*,
    mark size=2pt,
    color=OIblue,
] table {
x y
24.78 860
24.91 859
24.95 828
25.07 826
25.09 823
25.16 822
25.19 820
25.21 819
25.24 818
25.26 817
31.97 816
32.11 778
32.51 747
32.67 737
32.78 736
39.16 725
39.3 723
39.55 721
39.97 713
46.1 710
46.43 708
46.56 702
56.48 691
61.85 687
62.07 685
68.13 668
68.16 667
68.21 666
68.25 665
116.29 664
116.53 662
};
\addlegendentry{Baseline}

\addplot[
    const plot mark left,
    mark=square*,
    mark size=2pt,
    color=OIvermillion,
] table {
x y
13.64 724
13.66 723
13.69 718
23.89 715
24.09 710
24.11 708
24.88 704
24.91 698
24.93 678
33.4 676
33.46 671
47.85 670
83.08 668
84.07 666
84.1 665
84.44 662
118.24 661
};
\addlegendentry{PDHG}
\end{axis}
\end{tikzpicture}
\end{subfigure}

\vspace{0.8em}

\begin{subfigure}{0.45\textwidth}
\centering
\begin{tikzpicture}
\begin{axis}[
    grid=both,
    legend style={at={(0.98,0.98)},anchor=north east},
    mark size=1pt,
    xlabel={Runtime (s)},
    ylabel={Incumbent objective},
    ylabel style={yshift=5pt},
    legend image post style={mark size=2pt},
    title={rail03},
    ]

\definecolor{OIblue}{RGB}{0,114,178}
\definecolor{OIvermillion}{RGB}{213,94,0}

\addplot[
    const plot mark left,
    mark=*,
    mark size=2pt,
    color=OIblue,
] table {
x y
529.72 -201.955
572.58 -202.039
572.76 -485.783
577.98 -512.205
581.88 -525.122
587.59 -557.783
592.69 -831.622
600.92 -835.05
606.13 -845.127
629.08 -845.244
631.84 -848.561
642.16 -848.594
642.18 -848.627
650.25 -848.661
657.88 -848.827
723.67 -851.65
723.69 -854.661
1222.55 -854.777
};
\addlegendentry{Baseline}

\addplot[
    const plot mark left,
    mark=square*,
    mark size=2pt,
    color=OIvermillion,
] table {
x y
268.21 -763.788
271.36 -770.094
271.38 -772.866
271.4 -772.883
271.43 -773.133
272.87 -804.916
325.17 -804.938
359.53 -805.022
366.74 -811.422
367.39 -811.45
367.94 -818.833
371.14 -822.111
375.91 -822.294
377.79 -824.372
377.82 -831.994
377.84 -832.016
377.86 -837.038
377.88 -837.238
404.36 -840.777
442.76 -850.594
577.29 -850.683
};
\addlegendentry{PDHG}
\end{axis}
\end{tikzpicture}
\end{subfigure}
\hfill
\begin{subfigure}{0.45\textwidth}
\centering
\begin{tikzpicture}
\begin{axis}[
    grid=both,
    legend style={at={(0.98,0.98)},anchor=north east},
    mark size=1pt,
    xlabel={Runtime (s)},
    ylabel style={yshift=3pt},
    legend image post style={mark size=2pt},
    title={ramos3},
    ]

\definecolor{OIblue}{RGB}{0,114,178}
\definecolor{OIvermillion}{RGB}{213,94,0}

\addplot[
    const plot mark left,
    mark=*,
    mark size=2pt,
    color=OIblue,
] table {
x y
0 542
6.41 437
6.41 310
6.41 309
6.42 300
6.42 296
6.42 295
9.51 294
9.51 293
9.62 292
9.62 290
19.51 289
19.51 271
19.51 269
19.52 265
19.62 263
19.84 262
19.84 261
19.84 258
19.84 257
19.88 253
20.39 252
20.39 251
20.48 250
20.48 248
20.98 247
};
\addlegendentry{Baseline}

\addplot[
    const plot mark left,
    mark=square*,
    mark size=2pt,
    color=OIvermillion,
] table {
x y
0 542
0.03 440
0.03 312
0.03 304
0.03 297
0.03 296
0.03 295
0.03 294
0.03 292
0.15 290
0.15 289
0.15 286
0.15 284
0.16 283
0.16 282
0.27 276
0.27 275
0.36 274
};
\addlegendentry{PDHG}
\end{axis}
\end{tikzpicture}
\end{subfigure}

\caption{Incumbent solution progress (everything was run on a 64-core AMD EPYC 9575F CPU, except for PDHG which was run on an RTX PRO 6000 GPU).}
\label{fig:anecdotal_miplib2}
\end{figure}
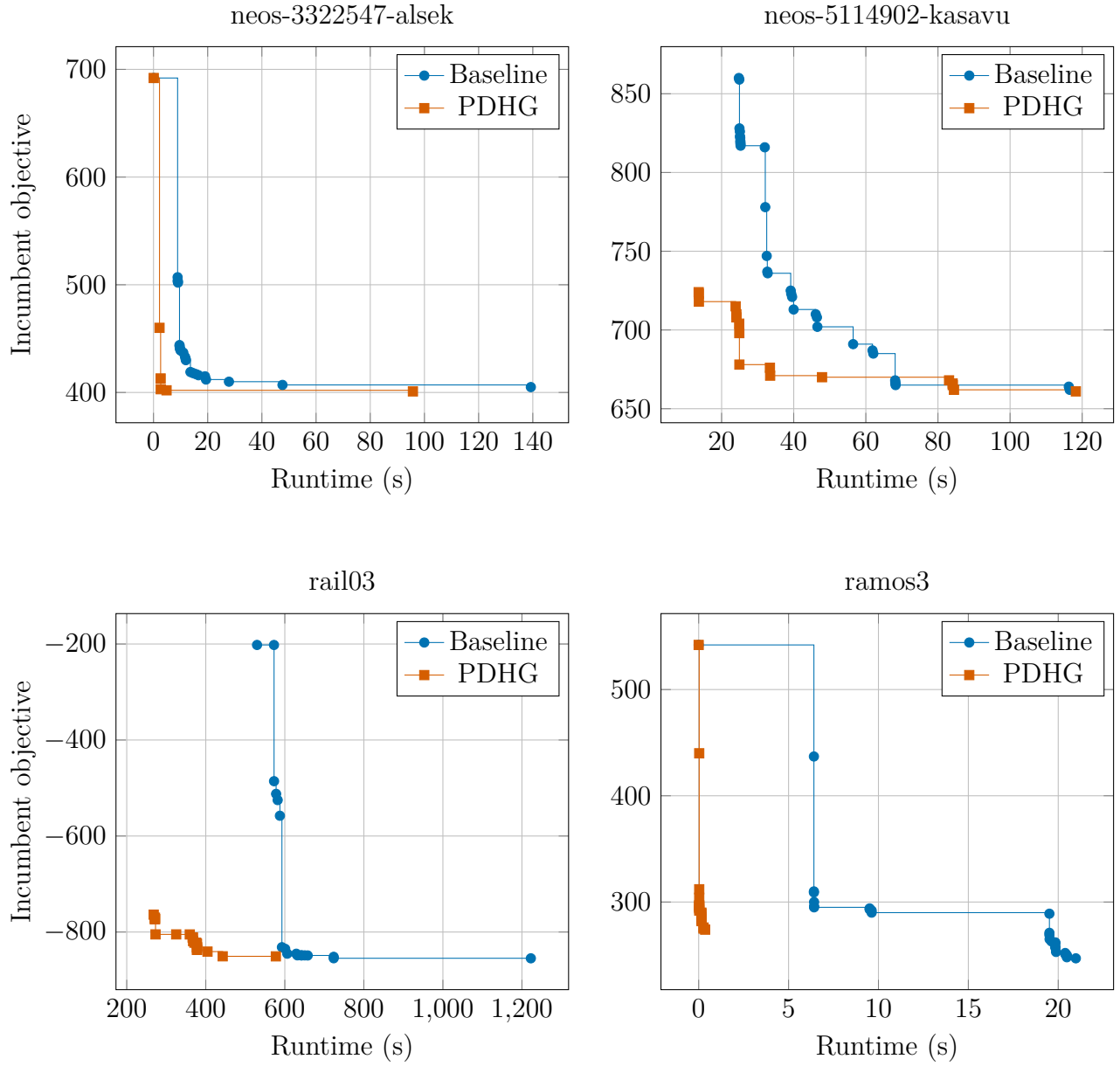

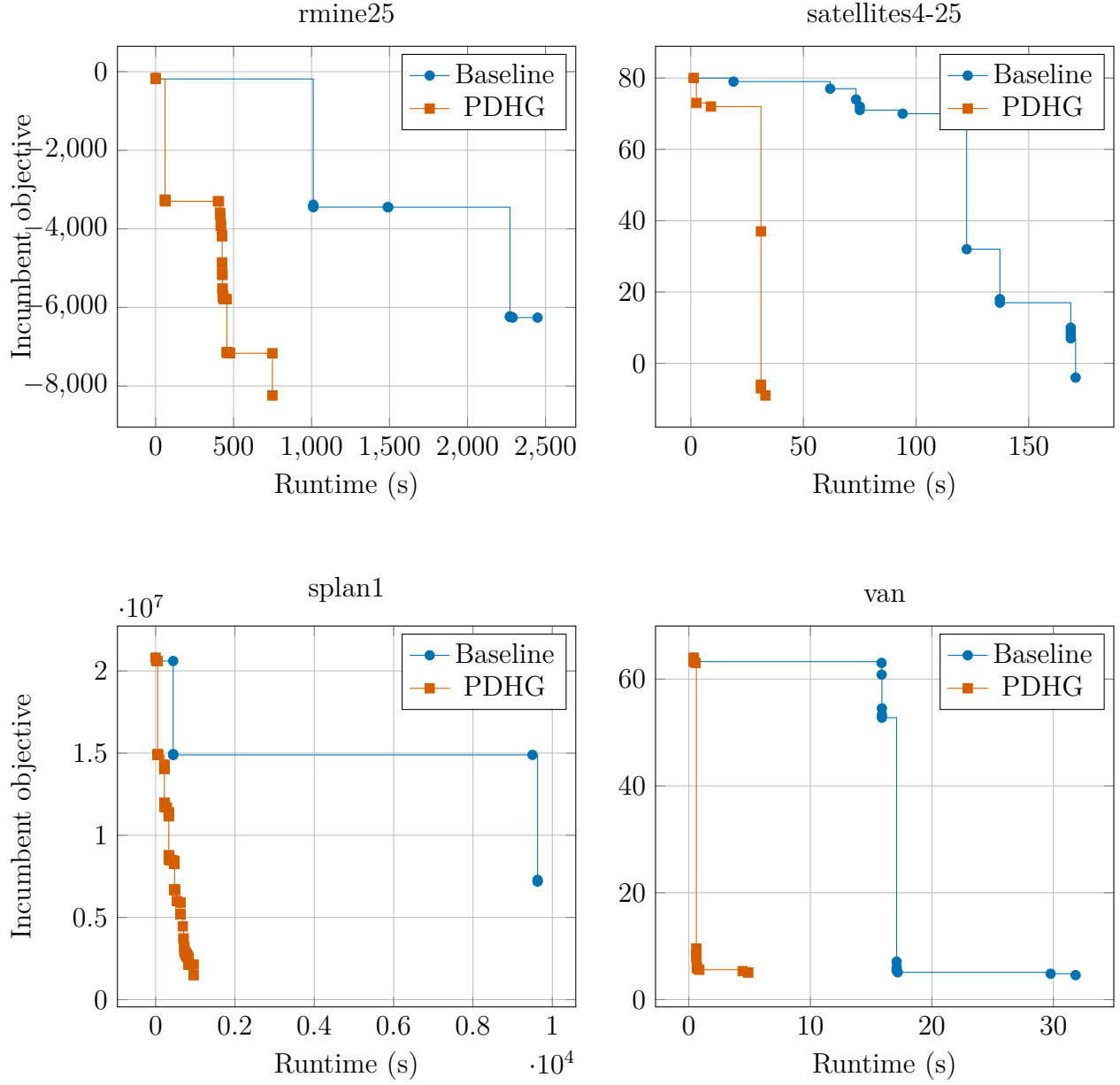
\begin{figure}[htbp]
\centering
\begin{subfigure}{0.45\textwidth}
\centering
\begin{tikzpicture}
\begin{axis}[
    grid=both,
    legend style={at={(0.98,0.98)},anchor=north east},
    mark size=1pt,
    xlabel={Runtime (s)},
    ylabel={Incumbent objective},
    ylabel style={yshift=5pt},
    legend image post style={mark size=2pt},
    title={rmine25},
    ]

\definecolor{OIblue}{RGB}{0,114,178}
\definecolor{OIvermillion}{RGB}{213,94,0}

\addplot[
    const plot mark left,
    mark=*,
    mark size=2pt,
    color=OIblue,
] table {
x y
0.09 -161.494
0.25 -185.328
1010.92 -3384.86
1011.04 -3435.41
1011.09 -3446.45
1487.82 -3448.44
1492 -3449.21
2272 -6228.43
2273 -6240.45
2273.06 -6240.47
2274.16 -6241.36
2276.83 -6241.42
2279.96 -6241.43
2281.64 -6241.64
2289.34 -6258.17
2449.17 -6259.59
};
\addlegendentry{Baseline}

\addplot[
    const plot mark left,
    mark=square*,
    mark size=2pt,
    color=OIvermillion,
] table {
x y
0.09 -161.494
0.25 -185.328
60.55 -3257.98
60.69 -3261.6
60.77 -3261.6
60.96 -3261.61
62.84 -3261.99
63.19 -3297.27
63.25 -3298.23
399.68 -3298.32
405.84 -3301
413.46 -3589.26
413.51 -3629.93
413.55 -3631.62
418.39 -3877.39
418.44 -3913.24
418.49 -3917.42
425.64 -4152.34
425.7 -4184.23
426.45 -4184.46
426.49 -4855.43
426.7 -4877.88
427.33 -5136.13
427.61 -5168.64
428.46 -5171.87
428.5 -5517.26
429.07 -5538.41
430.51 -5687.79
430.56 -5697.45
430.6 -5697.71
430.66 -5699.17
430.71 -5700
430.75 -5701.09
431.19 -5725.48
432.9 -5731.51
432.94 -5734.83
432.99 -5737.22
433.6 -5757.51
438.53 -5789.26
439.64 -5790.2
455.95 -5790.25
456.52 -7133.03
456.71 -7139.13
456.76 -7140.2
459.23 -7151.85
459.3 -7151.87
463.14 -7152.06
463.2 -7152.16
474.63 -7152.26
477.17 -7165.72
477.23 -7167.21
748.9 -7167.28
748.93 -8239.58
};
\addlegendentry{PDHG}
\end{axis}
\end{tikzpicture}
\end{subfigure}
\hfill
\begin{subfigure}{0.45\textwidth}
\centering
\begin{tikzpicture}
\begin{axis}[
    grid=both,
    legend style={at={(0.98,0.98)},anchor=north east},
    mark size=1pt,
    xlabel={Runtime (s)},
    ylabel style={yshift=5pt},
    legend image post style={mark size=2pt},
    title={satellites4-25},
    ]

\definecolor{OIblue}{RGB}{0,114,178}
\definecolor{OIvermillion}{RGB}{213,94,0}

\addplot[
    const plot mark left,
    mark=*,
    mark size=2pt,
    color=OIblue,
] table {
x y
1.29 80
18.97 79
61.87 77
73.35 74
75.03 72
75.03 71
94.05 70
122.45 32
137.21 18
137.22 17
168.65 10
168.65 9
168.66 8
168.66 7
170.82 -4
};
\addlegendentry{Baseline}

\addplot[
    const plot mark left,
    mark=square*,
    mark size=2pt,
    color=OIvermillion,
] table {
x y
1.31 80
2.49 73
8.94 72
31.15 37
31.16 -6
31.16 -7
33.14 -9
};
\addlegendentry{PDHG}
\end{axis}
\end{tikzpicture}
\end{subfigure}

\vspace{0.8em}

\begin{subfigure}{0.45\textwidth}
\centering
\begin{tikzpicture}
\begin{axis}[
    grid=both,
    legend style={at={(0.98,0.98)},anchor=north east},
    mark size=1pt,
    xlabel={Runtime (s)},
    ylabel={Incumbent objective},
    ylabel style={yshift=5pt},
    legend image post style={mark size=2pt},
    title={splan1},
    ]

\definecolor{OIblue}{RGB}{0,114,178}
\definecolor{OIvermillion}{RGB}{213,94,0}

\addplot[
    const plot mark left,
    mark=*,
    mark size=2pt,
    color=OIblue,
] table {
x y
0.13 2.0804e+07
16.93 2.06032e+07
442.01 2.06032e+07
442.1 1.49323e+07
442.15 1.49321e+07
442.91 1.49021e+07
443.01 1.48921e+07
9501.05 1.48921e+07
9629.14 7.28538e+06
9629.25 7.19534e+06
9629.29 7.19533e+06
};
\addlegendentry{Baseline}

\addplot[
    const plot mark left,
    mark=square*,
    mark size=2pt,
    color=OIvermillion,
] table {
x y
0.13 2.0804e+07
16.94 2.06032e+07
50.78 2.06032e+07
50.88 1.49323e+07
50.92 1.49321e+07
51.77 1.48939e+07
51.87 1.48839e+07
64.93 1.48839e+07
219.21 1.429e+07
219.68 1.42799e+07
219.76 1.41095e+07
219.82 1.40396e+07
219.91 1.40375e+07
220.8 1.19657e+07
220.86 1.19657e+07
221.18 1.19555e+07
221.83 1.18141e+07
222.09 1.1784e+07
224 1.1784e+07
224.06 1.1784e+07
224.11 1.1784e+07
224.16 1.1784e+07
224.22 1.1784e+07
224.27 1.1774e+07
224.33 1.17739e+07
225.07 1.17239e+07
284.44 1.16979e+07
330.36 1.1381e+07
331.2 1.1381e+07
331.31 1.13307e+07
331.38 1.11604e+07
331.48 1.11592e+07
331.97 8.77367e+06
332.03 8.77367e+06
332.46 8.77367e+06
333.2 8.73337e+06
333.26 8.73336e+06
333.7 8.68281e+06
333.79 8.67282e+06
333.93 8.67279e+06
333.99 8.67279e+06
334.51 8.63255e+06
334.57 8.63255e+06
335.4 8.62259e+06
335.46 8.62258e+06
336.43 8.57205e+06
336.76 8.54195e+06
337.51 8.54193e+06
337.97 8.54192e+06
339.81 8.54191e+06
339.87 8.54191e+06
339.94 8.54191e+06
345.36 8.54181e+06
345.42 8.54178e+06
345.49 8.54176e+06
345.9 8.49175e+06
381.63 8.49175e+06
410.44 8.47821e+06
410.56 8.46817e+06
465.47 8.43717e+06
468.02 8.42713e+06
468.32 8.41713e+06
468.39 8.30674e+06
468.53 8.30629e+06
468.6 8.30628e+06
473.13 8.25553e+06
475.5 8.2553e+06
478.14 6.70373e+06
478.22 6.70372e+06
479.96 6.70354e+06
480.02 6.70353e+06
480.09 6.67347e+06
480.51 6.6734e+06
480.57 6.67337e+06
480.64 6.67334e+06
535.4 5.99744e+06
624.97 5.91673e+06
625.17 5.90657e+06
625.24 5.21345e+06
625.3 5.21344e+06
625.39 5.20754e+06
625.45 5.19747e+06
684.86 4.45796e+06
702.63 3.71392e+06
702.69 3.71392e+06
721.54 3.18244e+06
721.61 3.17243e+06
725.31 2.97858e+06
737.95 2.86686e+06
752.06 2.86636e+06
752.13 2.85607e+06
757.29 2.85607e+06
757.35 2.77501e+06
767.58 2.77501e+06
767.64 2.77465e+06
767.71 2.73349e+06
777.01 2.72223e+06
777.08 2.72217e+06
777.14 2.71201e+06
777.21 2.70183e+06
777.27 2.70129e+06
777.34 2.69148e+06
777.42 2.68115e+06
777.49 2.68115e+06
787.97 2.68109e+06
788.04 2.68071e+06
788.1 2.67098e+06
788.17 2.6705e+06
801.29 2.65016e+06
801.36 2.64003e+06
814.53 2.58753e+06
814.6 2.55935e+06
814.67 2.55802e+06
825.87 2.55801e+06
837.84 2.16705e+06
837.96 2.16705e+06
838.85 2.1475e+06
839.04 2.14719e+06
839.14 2.14719e+06
839.87 2.14707e+06
840.04 2.14689e+06
840.6 2.13685e+06
840.67 2.13685e+06
842.45 2.13635e+06
842.8 2.1363e+06
843.86 2.13629e+06
843.99 2.13627e+06
958.33 2.13627e+06
958.38 1.48808e+06
};
\addlegendentry{PDHG}
\end{axis}
\end{tikzpicture}
\end{subfigure}
\hfill
\begin{subfigure}{0.45\textwidth}
\centering
\begin{tikzpicture}
\begin{axis}[
    grid=both,
    legend style={at={(0.98,0.98)},anchor=north east},
    mark size=1pt,
    xlabel={Runtime (s)},
    ylabel style={yshift=5pt},
    legend image post style={mark size=2pt},
    title={van}
    ]

\definecolor{OIblue}{RGB}{0,114,178}
\definecolor{OIvermillion}{RGB}{213,94,0}

\addplot[
    const plot mark left,
    mark=*,
    mark size=2pt,
    color=OIblue,
] table {
x y
0.4 64
0.41 63.5146
0.41 63.2648
15.86 63.0292
15.87 60.8378
15.89 54.5134
15.89 53.2928
15.89 52.7782
17.1 7.10262
17.1 6.08845
17.1 5.60305
17.11 5.35326
17.19 5.10348
29.77 4.83866
31.83 4.57385
};
\addlegendentry{Baseline}

\addplot[
    const plot mark left,
    mark=square*,
    mark size=2pt,
    color=OIvermillion,
] table {
x y
0.4 64
0.41 63.5146
0.41 63.2648
0.56 63.0292
0.62 9.56106
0.62 8.237
0.62 7.97219
0.65 6.38331
0.68 5.85369
0.86 5.58888
4.43 5.32406
4.89 5.05925
};
\addlegendentry{PDHG}
\end{axis}
\end{tikzpicture}
\end{subfigure}

\caption{Incumbent solution progress (everything was run on a 64-core AMD EPYC 9575F CPU, except for PDHG which was run on an RTX PRO 6000 GPU).}
\label{fig:anecdotal_miplib3}
\end{figure}

Figures~\ref{fig:anecdotal_miplib1},
\ref{fig:anecdotal_miplib2}, and
\ref{fig:anecdotal_miplib3}
show incumbent objective value versus time
for twelve MIPLIB models where our PDHG-based approach finds feasible
solutions much more quickly than the baseline.
\begin{figure}[htbp]
\centering
\begin{subfigure}{0.45\textwidth}
\centering
\begin{tikzpicture}
\begin{axis}[
    grid=both,
    legend style={at={(0.98,0.98)},anchor=north east},
    mark size=1pt,
    xlabel={Runtime (s)},
    ylabel={Incumbent objective},
    ylabel style={yshift=3pt},
    legend image post style={mark size=2pt},
    title={cdc7-4-3-2},
    ]

\definecolor{OIblue}{RGB}{0,114,178}
\definecolor{OIvermillion}{RGB}{213,94,0}

\addplot[
    const plot mark left,
    mark=*,
    mark size=2pt,
    color=OIblue,
] table {
x y
0 -232
119.6 -233
120.03 -235
120.03 -237
120.14 -238
120.14 -239
120.14 -240
120.14 -241
120.33 -243
121.26 -244
121.35 -246
121.45 -247
122.53 -248
122.72 -249
270.64 -250
270.64 -251
270.64 -252
};
\addlegendentry{Baseline}

\addplot[
    const plot mark left,
    mark=square*,
    mark size=2pt,
    color=OIvermillion,
] table {
x y
0 -232
0.26 -233
0.46 -239
0.46 -240
0.57 -241
0.57 -242
0.62 -246
0.62 -247
0.87 -249
0.89 -251
};
\addlegendentry{PDHG}
\end{axis}
\end{tikzpicture}
\end{subfigure}
\hfill
\begin{subfigure}{0.45\textwidth}
\centering
\begin{tikzpicture}
\begin{axis}[
    grid=both,
    legend style={at={(0.98,0.98)},anchor=north east},
    mark size=1pt,
    xlabel={Runtime (s)},
    ylabel style={yshift=5pt},
    legend image post style={mark size=2pt},
    title={graph20-80-1rand},
    ]

\definecolor{OIblue}{RGB}{0,114,178}
\definecolor{OIvermillion}{RGB}{213,94,0}

\addplot[
    const plot mark left,
    mark=*,
    mark size=2pt,
    color=OIblue,
] table {
x y
0 0
4.51 -1
4.51 -2
25.45 -3
25.45 -4
48.59 -5
69.76 -6
};
\addlegendentry{Baseline}

\addplot[
    const plot mark left,
    mark=square*,
    mark size=2pt,
    color=OIvermillion,
] table {
x y
0 0
0.29 -1
1 -2
2 -3
2.21 -4
3.78 -5
3.78 -6
};
\addlegendentry{PDHG}
\end{axis}
\end{tikzpicture}
\end{subfigure}

\vspace{0.8em}

\begin{subfigure}{0.45\textwidth}
\centering
\begin{tikzpicture}
\begin{axis}[
    grid=both,
    legend style={at={(0.98,0.98)},anchor=north east},
    mark size=1pt,
    xlabel={Runtime (s)},
    ylabel={Incumbent objective},
    ylabel style={yshift=5pt},
    legend image post style={mark size=2pt},
    title={graphdraw-grafo2},
    ]

\definecolor{OIblue}{RGB}{0,114,178}
\definecolor{OIvermillion}{RGB}{213,94,0}

\addplot[
    const plot mark left,
    mark=*,
    mark size=2pt,
    color=OIblue,
] table {
x y
13.29 1.48132e+06
145.34 963074
145.47 957674
145.48 956816
255.11 345662
255.11 344606
255.12 340584
354.03 264284
354.04 224776
354.04 224288
354.04 207738
354.05 201488
354.05 198434
354.06 197824
354.06 197820
354.06 197816
354.07 196832
354.07 195006
354.07 194528
354.08 194528
354.08 194520
354.08 192814
354.09 191824
354.09 189116
354.09 182512
354.1 182020
354.1 177876
354.11 175644
354.11 169078
354.11 167384
354.12 163946
354.12 163928
354.12 163920
354.13 163916
354.13 163420
354.13 163174
354.14 162432
354.14 158608
354.14 158486
354.15 158366
354.15 157892
354.15 156042
354.16 155918
354.16 155798
354.16 155684
354.17 154454
354.17 154448
354.18 153956
354.18 153588
354.18 153464
354.19 152384
354.19 151272
354.19 150410
465.03 147974
465.03 147638
465.03 147618
465.04 147616
465.04 147168
465.04 146790
465.05 146518
465.05 146516
465.05 146396
465.06 146394
465.06 146272
465.07 146150
465.07 146148
465.07 146026
465.08 146024
465.08 145536
465.08 145522
626.84 144936
626.84 144814
626.84 144778
626.85 144330
626.85 144322
626.85 144294
626.86 144286
626.86 143332
626.87 143328
626.87 143082
626.87 142712
626.88 142708
626.88 142626
626.88 141036
626.89 141032
626.89 140904
626.89 140898
626.9 140690
626.9 140682
768.95 140676
768.95 139814
768.96 138960
768.96 138220
768.96 137244
768.97 137240
768.97 136504
768.97 136496
768.98 136494
768.98 136400
963.59 133222
963.6 132980
963.6 132970
963.6 132142
963.61 132116
1227.28 131762
1227.28 131636
2285.47 131270
};
\addlegendentry{Baseline}

\addplot[
    const plot mark left,
    mark=square*,
    mark size=2pt,
    color=OIvermillion,
] table {
x y
13.29 1.48132e+06
21.67 188602
21.77 188102
21.78 187982
28.86 187980
28.86 187578
37.96 187570
37.96 184012
37.97 183282
37.97 182918
44.57 182792
44.57 178540
48.58 178534
48.59 175982
48.59 171562
48.59 171068
48.59 170704
48.6 170268
48.6 170144
48.6 170024
48.61 169902
48.61 169892
48.62 169768
51.56 168906
51.56 166580
54.72 163970
54.73 159942
54.73 159820
54.73 159818
54.74 159320
54.75 157856
54.8 157728
54.82 157722
54.82 157600
54.84 156484
54.84 156360
};
\addlegendentry{PDHG}
\end{axis}
\end{tikzpicture}
\end{subfigure}
\hfill
\begin{subfigure}{0.45\textwidth}
\centering
\begin{tikzpicture}
\begin{axis}[
    grid=both,
    legend style={at={(0.98,0.98)},anchor=north east},
    mark size=1pt,
    xlabel={Runtime (s)},
    ylabel style={yshift=5pt},
    legend image post style={mark size=2pt},
    title={in},
    ]

\definecolor{OIblue}{RGB}{0,114,178}
\definecolor{OIvermillion}{RGB}{213,94,0}

\addplot[
    const plot mark left,
    mark=*,
    mark size=2pt,
    color=OIblue,
] table {
x y
0.19 1489
14.92 1486
316.03 272
577.87 63
};
\addlegendentry{Baseline}

\addplot[
    const plot mark left,
    mark=square*,
    mark size=2pt,
    color=OIvermillion,
] table {
x y
0.19 1489
14.52 1486
17.99 64
};
\addlegendentry{PDHG}
\end{axis}
\end{tikzpicture}
\end{subfigure}

\caption{Incumbent solution progress (everything was run on a 64-core AMD EPYC 9575F CPU, except for PDHG which was run on an RTX PRO 6000 GPU).}
\label{fig:anecdotal_miplib1}
\end{figure}
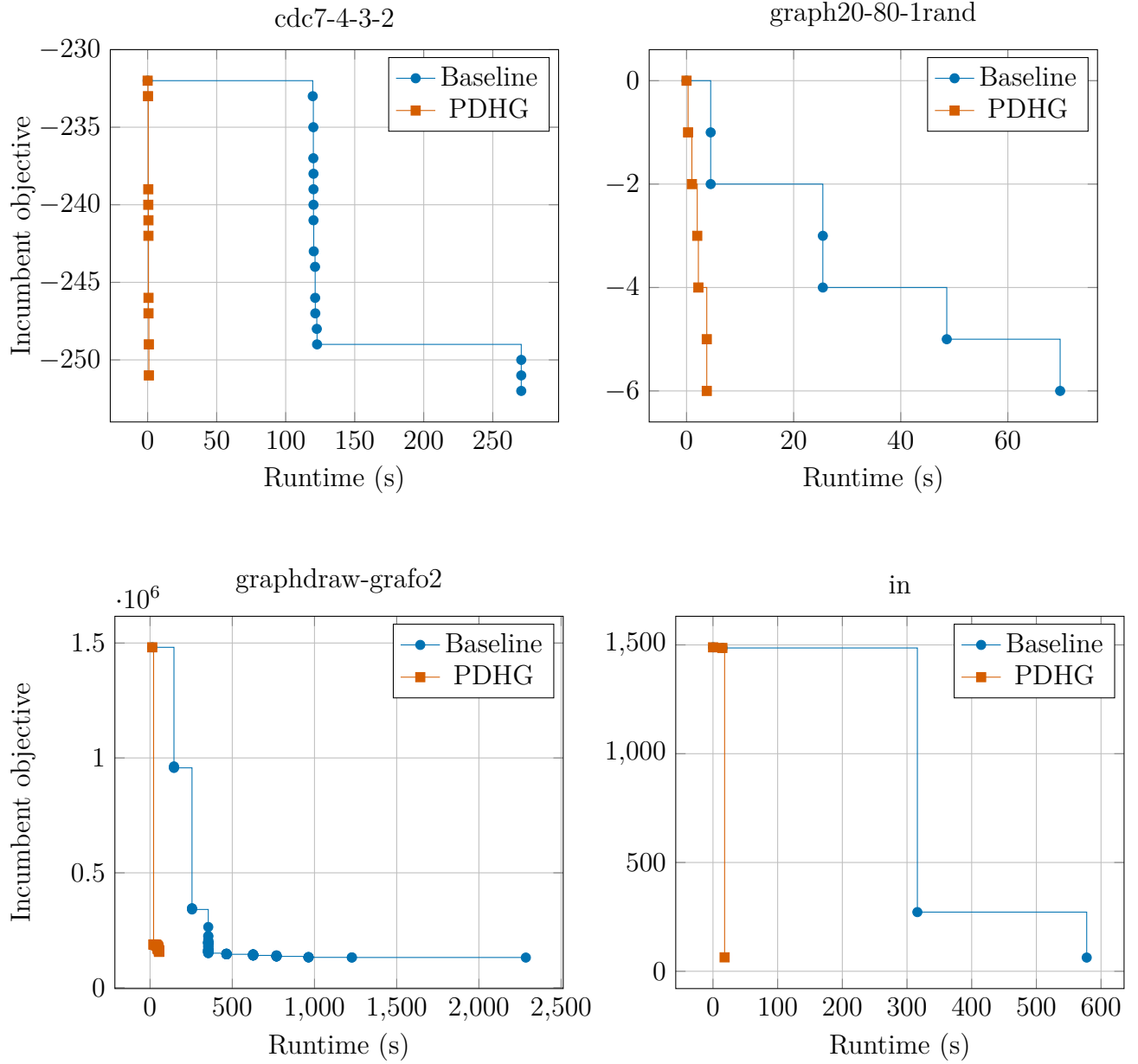

We actually saw improvements for roughly 20 of the 87
models in our test, but we expect that the success rate would be much
higher on a set with more difficult LP relaxations.  We find it quite
promising that this heuristic is able to frequently beat our
state-of-the-art baseline approach.

We have also seen dramatic improvements on customer models.  One
particular airline customer is solving a set of large,
extremely difficult MIP models.  Our PDHG-based heuristic approach
consistently finds near-optimal solutions 10 times faster than any of
our default options, and 2--10 times faster than an approach that uses
PDHG plus crossover at the root node.

One thing we would like to note about these results is that solutions
for the PDHG-based approach are found by a variety of Gurobi
heuristics.  In the case of the customer model, the best solution is
found with a very obscure heuristic.  This reinforces our point that,
if we had built a monolithic PDHG-based heuristic, rather than
integrating PDHG within nearly all of the existing Gurobi heuristics,
it would likely not have found many of these incumbents.

\section{Discussion}

One natural question at this point is whether it would be possible to
improve the results of our PDHG-based heuristic by finding ways to
reenable features that had to be turned off.  The natural candidate is
performing multiple root node solves in parallel to exploit
performance variability.  While it seems unlikely that we could run
multiple PDHG solves in parallel, given their heavy computing
requirements, it may be possible to capture intermediate PDHG
solutions and use them to launch heuristics on multiple threads.
Earlier PDHG iterates would be less accurate than the final, converged
iterate, but the results in this paper indicate
that low solution accuracy does not appear to significantly degrade
the quality of the MIP solutions that are produced.

Another option in a GPU context would be to launch a single PDHG-based
root cut thread within the baseline MIP solve.  We have found
that running PDHG on the GPU does not significantly degrade the
performance of tasks running on other CPU threads.  This option
could potentially exploit a different type of performance variability.
The software engineering required to manage both basic and non-basic
relaxations, using CPU cores and a GPU, in the same process would
be quite involved, but the ability to seamlessly integrate the
results from PDHG into a traditional MIP solve may be worth the
effort.

One downside of expensive MIP heuristics in general is that they focus
on only one side of the optimality gap: the incumbent objective.  Our
PDHG-based heuristic may prove useful for improving the lower bound as
well.  Recall that the cutting planes that are discovered at the PDHG
root node are valid for the original model (provided that we avoid invalid
reduced-cost fixings).  These cuts are often computed much more
quickly than with the default MIP solver.  We could transfer these
cuts to the default MIP, solve a single relaxation that includes them,
and continue the root solve from that point.

We should point out that there may be other ways to beat baseline
Gurobi on the models where our PDHG-based heuristic approach wins.
One option when crossover is the main expense would be to use an
interior-point solver without crossover.  Another is the Gurobi {\em
  NoRel\/} heuristic, which does not solve relaxations at all.
Another possibility is to integrate PDHG solutions into some of these
other schemes.  This will require further investigation.

\section{Conclusions}

This paper posed a simple question: what would happen if you replaced
the method used to solve LP relaxations within a state-of-the-art MIP
solver with PDHG?  While some features become degraded or disabled, we
have found that the resulting method retains most of Gurobi's ability
to combine a variety of powerful techniques to find high-quality
solutions.  We have shown that, for models where LP relaxations are
expensive to solve and where PDHG provides a substantial runtime
advantage, even low-accuracy relaxation solutions can provide
substantial improvements in the time required to find high-quality
feasible MIP solutions.

\section*{Acknowledgements}

The author would like to thank Nils-Christian Kempke and Lennart Lahrs
for their thoughtful feedback on an earlier draft of this paper, which
substantially improved the clarity of the presentation of results.
The author would also like to thank Chris Maes for the reminder that
even very noisy dual solution vectors provide valid objective bounds
for many LP relaxations.

\bibliographystyle{plainnat}
\bibliography{references}

\end{document}